\documentclass[10pt]{amsart}
            
\usepackage{amsmath}               
\usepackage{amsfonts}              
\usepackage{amsthm}                
\usepackage[utf8]{inputenc}
\usepackage{leftidx}
\usepackage{tikz-cd}
\usepackage{mathtools}

\newtheorem{tw}{Theorem}[section]
\newtheorem{lem}[tw]{Lemma}
\newtheorem{defi}[tw]{Definition}
\newtheorem{prop}[tw]{Proposition}
\newtheorem{rem}[tw]{Remark}
\newtheorem{cor}[tw]{Corollary}

\title{A Spectral Sequence for locally free isometric Lie group actions}
\author{Paweł Raźny}
\address{Institute of Mathematics \\
	Faculty of Mathematics and Computer Science \\
	Jagiellonian University in Cracow
	}
\email{pawel.razny@uj.edu.pl}


\begin{document}

\begin{abstract} We present a spectral sequence for free isometric Lie algebra actions (and consequently locally free isometric Lie group actions) which relates the de Rham cohomology of the manifold with the Lie algebra cohomology and basic cohomology (intuitively the cohomology of the orbit space). In the process of developing this sequence we introduce a new description of the de Rham cohomology of manifolds with such actions which appears to be well suited to this and similar problems. Finally, we provide some simple applications generalizing the Wang long exact sequence to Lie algebra actions of low codimension.
\end{abstract}
\maketitle
\section{Introduction}
The main goal of this paper is to provide a spectral sequence for isometric Lie group actions connecting the de Rham cohomology of the manifolds with the Lie algebra cohomology (describing the cohomological behaviour in the direction tangent to the orbits) and the basic cohomology (i.e. the cohomology of the orbit space) of the foliation induced by the orbits (describing the cohomological behaviour in the direction transverse to the orbits). This is achieved by the following Theorem:
\begin{tw}\label{MainG} Let $(M^{n+s},g)$ be a compact manifold with an isometric locally free action of an $s$-dimensional connected Lie group $G$. Then, there is a spectral sequence $E^{p,q}_r$ with:
\begin{itemize}
\item $E^{p,q}_2=H^p(M\slash\mathcal{F}, H^{q}(\mathfrak{g}))$, where $\mathfrak{g}$ is the Lie algebra of $G$ and $\mathcal{F}$ is the foliation generated by the fundamental vector fields of the action $\xi_1,...,\xi_s$.
\item $E^{p,q}_r$ converges to $H^{\bullet}_{dR}(M)$.  
\end{itemize}
\end{tw}
The above theorem can be easily deduced from its version for Lie algebra actions by simply passing to the Lie algebra action corresponding to the given Lie group action. 
\begin{tw}\label{MainA} Let $(M^{n+s},g)$ be a compact manifold with a Killing free action of an $s$-dimensional Lie algebra $\mathfrak{g}$. Then, there is a spectral sequence $E^{p,q}_r$ with:
\begin{itemize}
\item $E^{p,q}_2=H^p(M\slash\mathcal{F}, H^{q}(\mathfrak{g}))$, where $\mathcal{F}$ is the foliation generated by the fundamental vector fields of the action $\xi_1,...,\xi_s$.
\item $E^{p,q}_r$ converges to $H^{\bullet}_{dR}(M)$.  
\end{itemize}
\end{tw}
These results can be viewed as a common generalization of the following results:
\begin{enumerate}
\item The special case of Leray-Serre spectral sequence for principal $G$-bundles for a compact connected Lie group $G$. Here the orbit space is a manifold and hence for a principal $G$-bundle $P$ the second page is simply $H^p_{dR}(P\slash G, H^{q}(\mathfrak{g}))$ while the sequence converges to $H^{\bullet}_{dR}(P)$.
\item The Lie algebra version of the Lyndon-Hochschild-Serre spectral sequence for $\mathfrak{h}$ a Lie algebra of a compact connected group $H$ and $\mathfrak{g}$ an ideal in $\mathfrak{h}$. In this case the second page is $H^p(\mathfrak{h}\slash \mathfrak{g}, H^{q}(\mathfrak{g}))$ while the sequence converges to $H^{\bullet}(\mathfrak{h})\cong H^{\bullet}_{dR}(H)$. Note that for this spectral sequence to work the subgroup $G\subset H$ corresponding to $\mathfrak{g}$ need not be compact and hence it can be applied to some cases not covered by the previous point.
\item The Gysin long exact sequence for flows of nowhere vanishing Killing vector fields which was thoroughly studied in the case of Sasakian manifolds (cf. \cite{BG}) and recently studied in the general case (cf. \cite{W}) and even generalized to $\mathbb{S}^3$-actions (cf. \cite{RPSA}). These sequences (under the additional assumption that the action is locally free in the case of \cite{RPSA}) can be derived from our spectral sequence in a standard way.
\item The spectral sequence for almost $\mathcal{K}$-structures used in \cite{My3} to generalize results from $K$-contact and Sasakian manifold to this setting. In particular, the topological invariance of basic Betti numbers (cf. Theorem 7.4.14 from \cite{BG} for the Sasakian case), invariance under deformations of basic Hodge numbers (cf. \cite{My2} for the Sasakian case), and the classification of harmonic forms in terms of primitive basic harmonic forms (cf. Theorem 7.4.13 from \cite{BG} for the Sasakian case).
\end{enumerate}
\indent We prove Theorem \ref{MainA} in several steps. Firstly, in section \ref{frame} we choose a new metric $g'$ on $M$ so that the action remains isometric and the metric admits orthonormal fundamental vector fields $\{\xi_1,...,\xi_s\}$. Then, in section \ref{Tor} we study the behaviour of the forms $\{\eta_1,...,\eta_s\}$ dual to $\{\xi_1,...,\xi_s\}$ under the exterior derivative (Proposition \ref{split}). This enables us to establish a new cochain complex $\Omega_{\mathfrak{g}}^{\bullet}(M)\subset\Omega^{\bullet}(M)$ for such actions which is spanned by forms $\eta_{i_1}...\eta_{i_q}\alpha$ with $\alpha$ basic. This in turn allows us to easily compute the second page of the spectral sequence arising from the restriction of the standard filtration associated to a foliation (cf. \cite{ss1,ss2,ss3,ss4}) to $\Omega_{\mathfrak{g}}^{\bullet}(M)$. Sections \ref{SSplit} and \ref{Vanish} finish the proof by establishing the following theorem, which we believe to be of independent interest.
\begin{tw}\label{MainC} Let $(M^{n+s},g)$ be a compact manifold with a Killing free action of an $s$-dimensional Lie algebra $\mathfrak{g}$. Then, under the above notation, the inclusion $\Omega_{\mathfrak{g}}^{\bullet}(M)\subset\Omega^{\bullet}(M)$ induces an isomorphism in cohomology.
\end{tw}
In the final section of the paper we give some applications of our spectral sequence to low codimension actions, arriving at long exact sequences generalizing the Wang long exact sequence.

\bigskip
{\Small{\textbf{Acknowledgements} The author would like to thank Professor Martin Saralegui-Aranguren for his insightful remarks and sharing details of his work on related topics.}}
\bigskip
\section{Preliminaries}
\subsection{Foliations}
We provide a quick review of transverse structures on foliations.
\begin{defi} A codimension $q$ foliation $\mathcal{F}$ on a smooth $n$-manifold $M$ is given by the following data:
\begin{itemize}
\item An open cover $\mathcal{U}:=\{U_i\}_{i\in I}$ of M.
\item A q-dimensional smooth manifold $T_0$.
\item For each $U_i\in\mathcal{U}$ a submersion $f_i: U_i\rightarrow T_0$ with connected fibers (these fibers are called plaques).
\item For all intersections $U_i\cap U_j\neq\emptyset$ a local diffeomorphism $\gamma_{ij}$ of $T_0$ such that $f_j=\gamma_{ij}\circ f_i$
\end{itemize}
The last condition ensures that plaques glue nicely to form a partition of $M$ consisting of submanifolds of $M$ of codimension $q$. This partition is called a foliation $\mathcal{F}$ of $M$ and the elements of this partition are called leaves of $\mathcal{F}$.
\end{defi}
We call $T=\coprod\limits_{U_i\in\mathcal{U}}f_i(U_i)$ the transverse manifold of $\mathcal{F}$. The local diffeomorphisms $\gamma_{ij}$ generate a pseudogroup $\Gamma$ of transformations on $T$ (called the holonomy pseudogroup). The space of leaves $M\slash\mathcal{F}$ of the foliation $\mathcal{F}$ can be identified with $T\slash\Gamma$.
\begin{defi}
 A smooth form $\omega$ on $M$ is called transverse if for any vector field $X\in\Gamma (T\mathcal{F})$ (where $T\mathcal{F}$ denotes the bundle tangent to the leaves of $\mathcal{F}$) it satisfies $\iota_X\omega=0$. Moreover, if additionaly $\iota_Xd\omega=0$ holds for all $X$ tangent to the leaves of $\mathcal{F}$ then $\omega$ is said to be basic. Basic $0$-forms will be called basic functions henceforth.
\end{defi}
Basic forms are in one to one correspondence with $\Gamma$-invariant smooth forms on $T$. It is clear that $d\omega$ is basic for any basic form $\omega$. Hence, the set of basic forms of $\mathcal{F}$ (denoted $\Omega^{\bullet}(M\slash\mathcal{F})$) is a subcomplex of the de Rham complex of $M$. We define the basic cohomology of $\mathcal{F}$ to be the cohomology of this subcomplex and denote it by $H^{\bullet}(M\slash\mathcal{F})$. A transverse structure to $\mathcal{F}$ is a $\Gamma$-invariant structure on $T$. Among such structures the following is most relevant to our work:
\begin{defi}
$\mathcal{F}$ is said to be Riemannian if $T$ has a $\Gamma$-invariant Riemannian metric. This is equivalent to the existence of a Riemannian metric $g$ (called the transverse Riemannian metric) on $N\mathcal{F}:=TM\slash T\mathcal{F}$ with $\mathcal{L}_Xg=0$ for all vector fields $X$ tangent to the leaves.
\end{defi}
This structure enables the construction of a transverse version of Hodge theory (see \cite{E1}). Firstly, we recall a special class of Riemannian foliations on which the aforementioned theory is greatly simplified:
\begin{defi} A codimension $q$ foliation $\mathcal{F}$ on a compact connected manifold $M$ is called homologically orientable if $H^{q}(M\slash\mathcal{F})=\mathbb{R}$. A foliation $\mathcal{F}$ on a manifold $M$ is called homologically orientable if its restriction to each connected component of $M$ is homologically orientable.
\end{defi}
We will see in the final section that there is a large class of foliations considered in this paper which are homologically orientable and hence we shall restrict our attention to this case throughout the rest of this subsection.
\newline\indent Let $\mathcal{F}$ be a homologically orientable Riemannian foliation on a compact manifold $M$. One can use the transverse Riemannian metric to define the basic Hodge star operator $*_b:\Omega^k(M\slash\mathcal{F})\to\Omega^{q-k}(M\slash\mathcal{F})$ pointwise. This in turn allows us to define the basic adjoint operator:
$$\delta_b=(-1)^{q(r+1)+1}*_bd*_b.$$
\begin{rem} While we choose this to be the definition of $\delta_b$, it is in fact an adjoint of $d$ with respect to an appropriate inner product on forms induced by the transverse metric $g$. However, the definition of this inner product is quite involved and not necessary for our purpose. Although, we shall state some of the classical results of basic Hodge theory which use this inner product in their proof. For details see \cite{E1}.
\end{rem}
Using $\delta_b$ we can define the basic Laplace operator via:
$$\Delta_b=d\delta_b+\delta_bd.$$
As it turns out this operator has some nice properties similar to that of the classical Laplace operator. In particular, it is transversely elliptic in the following sense:
\begin{defi}
A basic differential operator of order $m$ is a linear map $D:\Omega^{\bullet}(M\slash\mathcal{F})\rightarrow\Omega^{\bullet}(M\slash\mathcal{F})$ such that in local coordinates $(x_1,...,x_p,y_1,...,y_q)$ (where $x_i$ are leaf-wise coordinates and $y_j$ are transverse ones) it has the form:
\begin{equation*}
D=\sum\limits_{|s|\leq m}a_s(y)\frac{\partial^{|s|}}{\partial^{s_1}y_1...\partial^{s_q}y_q}
\end{equation*}
where $a_s$ are matrices of appropriate size with basic functions as coefficients. A basic differential operator is called transversely elliptic if its principal symbol is an isomorphism at all points of $x\in M$ and all non-zero, transverse, cotangent vectors at x.
\end{defi}
In particular, this implies the following important result from \cite{E1}:
\begin{tw} Let $\mathcal{F}$ be a Riemannian homologically orientable foliation on a compact manifold $M$. Then:
\begin{enumerate}
\item $H^{\bullet}(M\slash\mathcal{F})$ is isomorphic to the space of basic harmonic forms $Ker(\Delta_b)$. In particular, it is finitely dimensional.
\item The basic Hodge star induces an isomorphism between $H^{k}(M\slash\mathcal{F})$ and $H^{q-k}(M\slash\mathcal{F})$ given by taking the class of the image through $*_b$ of a harmonic representative.
\end{enumerate}
\end{tw}
We finish this section by recalling the spectral sequence of a Riemannian foliation.
\begin{defi}\label{Fili} We put:
$$F^k_{\mathcal{F}}\Omega^r(M):=\{\alpha\in\Omega^r(M)\text{ }|\text{ } \iota_{X_{r-k+1}}...\iota_{X_1}\alpha=0,\text{ for } X_1,...,X_{r-k+1}\in\Gamma(T\mathcal{F})\}.$$
An element of $F^k_{\mathcal{F}}\Omega^r(M)$ is called an $r$-differential form of filtration $k$.
\end{defi}
The definition above in fact gives a filtration of the de Rham complex. Hence, via known theory from homological algebra we can cosntruct a spectral sequence as follows:
\begin{enumerate}
\item The $0$-th page is given by $E_0^{p,q}=F_{\mathcal{F}}^{p}\Omega^{p+q}(M)\slash F_{\mathcal{F}}^{p+1}\Omega^{p+q}(M)$ and $d^{p,q}_0:E_0^{p,q}\to E^{p,q+1}_0$ is simply the morphism induced by $d$.
\item The $r$-th page is given inductively by:
$$E_{r}^{p,q}:=Ker(d^{p,q}_{r-1})\slash Im(d^{p-r,q+r-1}_{r-1})=\frac{\{\alpha\in F^{p}_{\mathcal{F}}\Omega^{p+q}(M)\text{ }|\text{ } d\alpha\in F_{\mathcal{F}}^{p+r}\Omega^{p+q+1}(M)\}}{F^{p+1}_{\mathcal{F}}\Omega^{p+q}(M)+d(F_{\mathcal{F}}^{p-r+1}\Omega^{p+q-1}(M))}$$
\item The $r$-th coboundary operator $d_r:E_r^{p,q}\to E^{p+r,q-r+1}_r$ is again just the map induced by $d$ (due to the description of the $r$-th page this has the target specified above and is well defined).
\end{enumerate}
Furthermore, since the filtration is bounded this spectral sequence converges and its final page is isomorphic to the cohomology of the cochain complex (in this case the de Rham cohomology of $M$).
\begin{rem} The above spectral sequence can be thought of as a generalization of the Leray-Serre spectral sequence in de Rham cohomology to arbitrary Riemannian foliations (as opposed to fiber bundles).
\end{rem}
\subsection{Lie algebras and Lie groups}
In this subsection we give a brief overview of properties of Lie groups, Lie algebras and their actions which will be useful in the rest of the article. Throghout this section let $G$ be an $s$-dimensional Lie group and let $\mathfrak{g}$ be the Lie algebra of left invariant vector fields on $G$ which can be classically identified with $T_eG$ (where $e\in G$ denotes the identity element). Choose a basis $\{\tilde{\xi}_1,...,\tilde{\xi}_s\}$ of $\mathfrak{g}$ and let $\{Z_1,...,Z_s\}$ be a set of right invariant vector fields on $G$ such that $(Z_i)_e=(\tilde{\xi}_i)_e$ then it is well known that:
\begin{enumerate}
\item $Z_i$ are the fundamental vector fields of the left action of $G$ on itself,
\item $[\tilde{\xi}_i,\tilde{\xi}_j]_e=-[Z_i,Z_j]_e$ for all $i,j\in \{1,...,s\}$,
\item $[\tilde{\xi}_i,Z_j]=0$ for all $i,j\in \{1,...,s\}$.
\end{enumerate}
In particular, this justifies the following convention:
\begin{defi} We say that an $\mathbb{R}$-linear map $A:\mathfrak{g}\to \Gamma(TM)$ is an action of the Lie algebra $\mathfrak{g}$ on a smooth manifold $M$ if $A([\tilde{\xi},\tilde{\zeta}])=-[A(\tilde{\xi}),A(\tilde{\zeta})]$ for all $\tilde{\xi},\tilde{\zeta}\in\mathfrak{g}$. Given a basis $\{\tilde{\xi}_1,...,\tilde{\xi}_s\}$ of $\mathfrak{g}$ we say that the vector fields $\{A(\tilde{\xi}_1),...,A(\tilde{\xi}_s)\}$ are the corresponding fundamental vector fields. We say that a Lie algebra action is free if $\{A(\tilde{\xi}_1),...,A(\tilde{\xi}_s)\}$ are linearly independent at each point. We say that a Lie algebra action is Killing (or isometric) if $\{A(\tilde{\xi}_1),...,A(\tilde{\xi}_s)\}$ are Killing vector fields.
\end{defi}
\begin{rem} A Lie group action of $G$ on $M$ induces a Lie algebra action of $\mathfrak{g}$ on $M$ with the same fundamental vector fields. This action is free if and only if the corresponding Lie group action is locally free (i.e. the isotropy groups of the action are discrete) and it is Killing if and only if $G$ acts by isometries.
\end{rem}
Assuming that $G$ is compact any orbit of a smooth action of $G$ on a manifold $M$ is equivarianlty diffeomorphic to $G\slash G_x$ where $G_x$ is the isotropy group at any point $x$ in the given orbit. Moreover, two orbits $Gx$ and $Gy$ are equivariantly diffeomorphic if and only if $G_x$ is conjugate to $G_y$ (in this case we say that the two orbits are of the same type; cf. \cite{B}). Taking $K_i$ to be the sum of all orbits of the given type we get a Whitney stratification of $M$ (in the sense of \cite{Str}) with a finite number of strata. Moreover, we can order the strata so that if $H_i\supset H_j$ are some representatives of the conjugacy classes of subgroups of $G$ corresponding to $K_i$ and $K_j$ then $i\geq j$. For later use we are going to need the following property of tubular neighbourhoods of strata:
\begin{prop}\label{tub} Given a strata $K_i$ of the action of a compact group $G$ on $M$ there exists an invariant tubular neighbourhood $U$ of $K_i$ contained in $N=\bigcup\limits^i_{j=1} K_j$ such that $N\backslash \overline{U}$ is equivariantly diffeomorphic to $N\backslash K_i$.
\begin{proof} This is a simple consequence of the existence of invariant tubular neighbourhoods applied to $K_i\subset N$.
\end{proof}
\end{prop}
We are also going to need the folowing version of the Slice Theorem (cf. \cite{B}):
\begin{tw} Let $G$ be a compact Lie group acting on a manifold $M$. For each orbit $Gx$ there exists a tubular neighbourhood $U$ equivariantly diffeomorphic to $G\times_{G_x} N_x(Gx)$ such that the image $S$ of $\{e\}\times_{G_x} N_x(Gx)$ through this diffeomorphism is a $G_x$-invariant manifold passing through $x$ (where we consider the normal space $N_x(Gx)$ of the orbit of $x$ at $x$ with the action of $G_x$ given by derivatives at $x$). In this case we say that $S$ is a slice through $x$ and that $U$ is a neighbourhood admitting a slice.
\end{tw}
\indent It is well known (cf. \cite{Hoch1}) that the automorphism group $Aut(G)$ (with compact open topology) of a connected Lie group is a Lie group itself. This group is relevant to our work due to the following result which summarizes the properties of dense subgroups in the case of interest to us.
\begin{tw}[Theorem 2.2.1 \cite{vE}]\label{Dense} Let $f:G\to L$ be an injective morphism of Lie groups such, that the image $f(G)$ is dense. If the inner automorphism group $I(G)$ is a closed subgroup of $Aut(G)$ then $L=\overline{f(G)}=\overline{f(Z(G))}\cdot f(G)$, where $Z(G)$ denotes the center of $G$.
\end{tw}
\begin{rem}
In this paper we usually treat the group $G$ as a subgroup of its closure and use the symbol $G$ for both the group with its subgroup and Lie group topology. Hence, for example in the above Theorem we drop $f$ from the assertion leaving $\overline{G}=\overline{Z(G)}\cdot G$. This abuse of notation simplifies the notation somewhat without creating much confusion and we will leave comments on the topology considered whenever required.
\end{rem}
Finally, we wish to recall the construction of Lie algebra cohomology. The cochain complex, $(\bigwedge\text{}^{\bullet}\mathfrak{g}^*,d_{\mathfrak{g}})$
where,
$$d_{\mathfrak{g}}\alpha(X_0,...,X_{n}):=\sum_{1\leq i<j\leq n}(-1)^{i+j}\alpha([X_i,X_j],X_0,...,\hat{X}_i,...,\hat{X}_j,...,X_n),$$
is called the Chevalley–Eilenberg complex and its cohomology $H^{\bullet}(\mathfrak{g})$ the Lie algebra cohomology of $\mathfrak{g}$. While this cohomology theory has found many applications throughout mathematics its most important property for our purpose is the following string of isomorphisms (for a compact connected Lie group $G$):
$$H^{\bullet}_{dR}(G)\cong H^{\bullet}_{G}(G)\cong H^{\bullet}(\mathfrak{g}),$$
where the middle term is the cohomology of the cochain complex $\Omega_G^{\bullet}(G)$ of left invariant forms on $G$. The first isomorphism is induced simply by the inclusion of forms and can be proven by finding a cochain map that induces the inverse mapping in cohomology (averaging a form over $G$). The second isomorphism follows simply from comparing $d_{\mathfrak{g}}$ with the exterior derivative $d$ restricted to left invariant forms.
\section{Free Killing Lie algebra actions}\label{frame}
Let us assume that $(M^{n+s},g)$ is a compact Riemannian manifold with a free Killing Lie algebra action of an $s$-dimensional Lie algebra $\mathfrak{g}$ with a choosen basis $\tilde{\xi}_1,...,\tilde{\xi}_s$ and corresponding fundamental vector fields $\xi_1,...,\xi_s$. In this section we will show how to construct a new metric $g'$ on $M$ so that the action of $\mathfrak{g}$ admits orthonormal fundamental fields $\{\xi'_1,...,\xi'_s\}$ and remains a Killing action.
\newline\indent Firstly, let us notice that the setup implies that:
\begin{enumerate}
\item The splitting $TM:=T\mathcal{F}\oplus T\mathcal{F}^{\perp}$ is preserved by the corresponding group action which in turn implies that it is preserved by the operation $[\xi_i,\bullet]$.
\item The above splitting induces a splitting of $g$ into the transverse part $g_{\perp}$ and the tangent part $g_T$.
\item The transverse part is well known to induce a transverse Riemannian structure for the corresponding foliation and hence $\mathcal{L}_{\xi_i}g_{\perp}=0$.
\item From the previous point we conclude that $\mathcal{L}_{\xi_i}g_T=\mathcal{L}_{\xi_i}(g_{\perp}+g_T)=\mathcal{L}_{\xi_i}g=0$.
\end{enumerate}
Since the action of $\mathfrak{g}$ is Killing the subgroup $G\subset Diff(M)$ generated by the flows of $\{\xi_1,...,\xi_s\}$ is contained in the group of isometries $Isom(M,g)$ which is well known to be a compact (finitely dimensional) Lie group and hence, the closure $\overline{G}$ of $G$ in $Diff(M)$ is compact as well. This in turn implies that there is a bi-invariant metric $\tilde{g}'$ on $G$ given by the restriction of a bi-invariant metric on $\overline{G}$ (which exists by compactness). Put $\tilde{g}'_{ij}:=\tilde{g}'(\tilde{\xi}_i,\tilde{\xi}_j)$. This induces a new metric $g'_T$ on $T\mathcal{F}$ defined by the formula $g'_T(\xi_i,\xi_j)=\tilde{g}'_{ij}$. Finally, we put $g':=g_{\perp}+g'_T$.
\begin{lem} The vector fields $\{\xi_1,...,\xi_s\}$ are Killing with respect to $g'$.
\begin{proof} We consider the tensor:
 $$(\mathcal{L}_{\xi_i}g')(X_1,X_2):=\mathcal{L}_{\xi_i}g'(X_1,X_2)-g'([\xi_i,X_1],X_2)-g'(X_1,[\xi_i,X_2]).$$
\indent Firstly, let us assume that $X_1,X_2\in T\mathcal{F}^{\perp}$. In this case the brackets in the expression above are also elements of $T\mathcal{F}^{\perp}$ and hence we see that $(\mathcal{L}_{\xi_i}g')(X_1,X_2)=(\mathcal{L}_{\xi_i}g_{\perp})(X_1,X_2)=0$.
\newline\indent Secondly, for $X_1,X_2\in T\mathcal{F}$ a similar argument shows that $(\mathcal{L}_{\xi_i}g')(X_1,X_2)=(\mathcal{L}_{\xi_i}g'_T)(X_1,X_2)=0$ (where the second equality follows since the given metric restricted to orbits is induced by the bi-invariant metric on $G$).
\newline\indent Finally, for $X_1\in T\mathcal{F}^{\perp}$ and $X_2\in T\mathcal{F}$ each term on the right hand side vanishes.
\end{proof}
\end{lem}
Finally, we find vector fields $\{\xi'_1,...,\xi'_s\}$ by taking the fundamental fields corresponding to a basis $\{\tilde{\xi}'_1,...,\tilde{\xi}'_s\}$ of $\mathfrak{g}$ orthogonal with respect to $\tilde{g}'$.
\begin{rem} Throughout the paper we will assume at various points that the action of $\mathfrak{g}$ admits an orthonormal frame as above. However, it is important to note that via applying apriori the above construction our spectral sequence can be constructed for an arbitrary free Killing Lie algebra action. 
\end{rem}
\section{The spectral sequence of a free Killing Lie algbera action}\label{Tor}
Throuoghout this section let $(M,g)$ be an $(n+s)$-dimensional compact Riemannian manifold with a free Killing action of a Lie algebra $\mathfrak{g}$ with orthonormal fundamental vector fields $\xi_1,...,\xi_s$. Let us note that the decomposition $TM=T\mathcal{F}\oplus T\mathcal{F}^{\perp}$ induces a decomposition of the cotangent bundle:
$$T^*M=Ann(T\mathcal{F})\oplus Ann(T\mathcal{F}^{\perp})\cong T^*\mathcal{F}^{\perp}\oplus T^*\mathcal{F},$$
where $Ann(E_1)\subset E_2^*$ denotes the annihilator of the vector subbundle $E_1\subset E_2$. This in turn induces a bigrading on differential forms:
$$\Omega^{r}(M)=\bigoplus\limits_{p+q=r} (\bigwedge\text{}^p Ann(T\mathcal{F})\otimes\bigwedge\text{}^q Ann(T\mathcal{F}^{\perp})).$$
Furthermore, note that $Ann(T\mathcal{F}^{\perp})$ is pointwise spanned by forms $\eta_i:=g(\xi_i,\bullet)$ which allows us to write a form $\alpha$ of bidegree $(p,q)$ as:
$$\alpha=\sum_{1\leq i_1\leq...\leq i_q\leq s}\eta_{i_1}...\eta_{i_q}\alpha_{i_1...i_q},$$
for some transverse $p$-forms $\alpha_{i_1...i_q}$. We now wish to study how the operator $d$ behaves with respect to this bigradation. To this end we use the well known formula:
$$(d\alpha)(X_0,...,X_{r})=\sum_{i} (-1)^i\mathcal{L}_{X_i}\alpha(X_0,...,\hat{X}_i,..,X_r)+\sum_{i<j}(-1)^{i+j}\alpha([X_i.X_j],X_0,...,\hat{X}_i,...,\hat{X}_j,...,X_r).$$
Using this formula we see that the operator $d$ splits into $4$ components:
$$d=d^{(-1,2)}+d^{(0,1)}+d^{(1,0)}+d^{(2,-1)}.$$
Moreover, due to involutivity of $T\mathcal{F}$ it can be easilly seen that $d^{(-1,2)}=0$. It remains to show how the three remaining components act on $\eta_i$:
\begin{prop}\label{split} Under the above notation:
\begin{enumerate}
\item The cochain complex $(\bigwedge^{\bullet} span(\eta_1,...,\eta_s),d^{(0,1)})$ is isomorphic to the Lie algebra cochain complex of $\mathfrak{g}$ (with $-d_{\mathfrak{g}}$).
\item $d^{(1,0)}\eta_i=0$ for each $i\in\{1,...,s\}$.
\item $d^{(2,-1)}\eta_i$ is basic for each $i\in\{1,...,s\}$.
\end{enumerate}
\begin{proof} For the first point, let $\tilde{\xi}_i$ denote the vector in $\mathfrak{g}$ that maps to $\xi_i$ through the Lie algebra action. Due to our setup, the Riemannian metric $g$ induces a scalar product $\tilde{g}$ on $\mathfrak{g}$ for which $\{\tilde{\xi}_1,...,\tilde{\xi}_s\}$ is an orthonormal basis. Put $\tilde{\eta}_i:=\tilde{g}(\tilde{\xi}_i,\bullet)$. Then it is easy to check that the algebra morphism $A$ defined by $A\eta_i=\tilde{\eta}_i$ induces the desired isomorphism of chain complexes.
\newline\indent For the second point, we choose a vector field $X\in\Gamma(T\mathcal{F}^{\perp})$ and compute for any $i\in\{1,...,s\}$:
$$(d\eta_j)(X,\xi_i)=\mathcal{L}_X\eta_j(\xi_i)-\mathcal{L}_{\xi_i}\eta_j(X)-\eta_{j}([X,\xi_i])=\eta_j(\mathcal{L}_{\xi_i}X)=0.$$
\indent For the third part, we shall work locally around a point $x_0$. We start by lifting the Lie algebra action to a Lie group action of the group $G$ generated by the flows of $\{\xi_1,...,\xi_s\}$. Next we choose a small open neigbourhood $U_0$ of $x_0$ and let $T\subset U_0$ be a disc transverse to the orbits and passing through $x_0$. We also choose a small neighbourhood $V$ of $e$ in $G$, such that $V\times T$ is diffeomorphic through the action to a small neighbourhood $U\subset U_0$ of $x_0$. We can now use this diffeomorphism to define vector fields $X_1,...,X_s$ by pushing forward the left-invariant vector fields restricted to $V$, while in such coordinates $\xi_1,...,\xi_s$ are pushforwards of right invariant vector fields, which implies $\mathcal{L}_{X_i}\xi_{j}=0$ for all $i,j\in\{1,...,s\}$. On the other hand, we get $\mathcal{L}_{X_i}g=0$, since $\mathcal{L}_{\xi_i}g_{\perp}=0$ (for all $i\in\{1,...,s\}$) implies $\mathcal{L}_{X_i}g_{\perp}=0$ (for all $i\in\{1,...,s\}$) and we have:
$$0=\mathcal{L}_{X_i}g(\xi_j,\xi_k)=g([X_j,\xi_j],\xi_k)+g(\xi_j,[X_i,\xi_k])+(\mathcal{L}_{X_i}g)(\xi_j,\xi_k)=(\mathcal{L}_{X_i}g)(\xi_j,\xi_k).$$ 
This means that $\mathcal{L}_{X_i}\eta_j=0$ for $i,j\in\{1,...,s\}$, since $\mathcal{L}_{X_i}g=0$ and $[X_i,\xi_j]=0$. This in turn implies that:
$$\mathcal{L}_{X_i}d\eta_j=d\mathcal{L}_{X_i}\eta_j=0.$$
However, due to the decomposition of forms and the second point we have:
$$d^{(2,-1)}\eta_j=d\eta_j-d^{(0,1)}\eta_j.$$
Using the first point and the above consideration we see that:
$$\mathcal{L}_{X_i}d^{(2,-1)}\eta_j=\mathcal{L}_{X_i}(d\eta_j-d^{(0,1)}\eta_j)=0.$$
This in turn coupled with the transversality of $d^{(2,-1)}\eta_j$ proves that it is in fact basic (when restricted to $U$). Since $x_0$ was arbitrary this ends the proof.
\end{proof}
\end{prop}
This has a useful corollary which is the first step to defining our spectral sequence.
\begin{cor} The graded vector subspace $\Omega^{\bullet}_{\mathfrak{g}}(M)\subset\Omega^{\bullet}(M)$ spanned by differential forms of the form $\eta_{i_1}...\eta_{i_k}\alpha$ where $\alpha$ is basic, is in fact a subcomplex of $\Omega^{\bullet}(M)$.
\begin{proof} This follows directly from the previous Proposition since it implies that in the expression:
$$d(\eta_{i_1}...\eta_{i_k}\alpha)=\sum_j (-1)^{j-1}(d^{(0,1)}\eta_{i_j}+d^{(2,-1)}\eta_{i_j})\eta_{i_1}...\hat{\eta}_{i_j}...\eta_{i_k}\alpha+(-1)^{k}\eta_{i_1}...\eta_{i_k}d\alpha,$$
every summand is an element of $\Omega^{\bullet}_{\mathfrak{g}}(M)$ whenever $\eta_{i_1}...\eta_{i_k}\alpha\in\Omega^{\bullet}_{\mathfrak{g}}(M)$ (by Proposition \ref{split}).
\end{proof} 
\end{cor}
With this we can now define a filtration on $\Omega^{\bullet}_{\mathfrak{g}}(M)$ by restricting the filtration from Definition \ref{Fili} to such forms (in a fashion similar to \cite{My3}). In other words we put:
$$F^k_{\mathcal{F}}\Omega^r_{\mathfrak{g}}(M):=\{\alpha\in\Omega^r_{\mathfrak{g}}(M)\text{ }|\text{ } \iota_{X_{r-k+1}}...\iota_{X_1}\alpha=0,\text{ for } X_1,...,X_{r-k+1}\in\Gamma(T\mathcal{F})\}.$$
In particular, we can again construct a spectral sequence as follows:
\begin{enumerate}
\item The $0$-th page is given by $E_0^{p,q}=F_{\mathcal{F}}^{p}\Omega^{p+q}_{\mathfrak{g}}(M)\slash F_{\mathcal{F}}^{p+1}\Omega^{p+q}_{\mathfrak{g}}(M)$ and $d^{p,q}_0:E_0^{p,q}\to E^{p,q+1}_0$ is simply the morphism induced by $d$.
\item The $r$-th page is given inductively by:
$$E_{r}^{p,q}:=Ker(d^{p,q}_{r-1})\slash Im(d^{p-r,q+r-1}_{r-1})=\frac{\{\alpha\in F^{p}_{\mathcal{F}}\Omega^{p+q}_{\mathfrak{g}}(M)\text{ }|\text{ } d\alpha\in F_{\mathcal{F}}^{p+r}\Omega^{p+q+1}_{\mathfrak{g}}(M)\}}{F^{p+1}_{\mathcal{F}}\Omega^{p+q}_{\mathfrak{g}}(M)+d(F_{\mathcal{F}}^{p-r+1}\Omega^{p+q-1}_{\mathfrak{g}}(M))}$$
\item The $r$-th coboundary operator $d_r:E_r^{p,q}\to E^{p+r,q-r+1}_r$ is again just the map induced by $d$.
\end{enumerate}
Furthermore, since the filtration is bounded this spectral sequence converges and its final page is isomorphic to the cohomology $H^{\bullet}_{\mathfrak{g}}(M)$ of the cochain complex $\Omega^{\bullet}_{\mathfrak{g}}(M)$. We shall show in the subsequent section that this is in fact isomorphic to the de Rham cohomology of $M$. We end this section by computing the second page of this spectral sequence.
\begin{prop} Under the above notation $E_2^{p,q}\cong H^p(M\slash\mathcal{F}, H^q(\mathfrak{g}))$.
\begin{proof} By $1)$ of Proposition \ref{split} the first page $E_1^{p,q}$ is equal to $\Omega^p(M\slash\mathcal{F})\otimes H^q(\mathfrak{g})$. Using $2)$ of Proposition \ref{split} we see that $d_1^{p,q}$ is simply the operator $d$ on basic forms with values in $H^q(\mathfrak{g})$.
\end{proof}
\end{prop}
\section{A splitting of de Rham cohomology}\label{SSplit}
Throughout this section let $(M^{n+s},g)$ be a Riemannian manifold on which there is a free Killing Lie algebra action. Without loss of generality (see Section \ref{frame}) we can assume that this action admits orthonormal fundamental fields $\{\xi_1,...,\xi_s\}$. Denote by $G$ the subgroup of $Isom(M,g)$ generated by the flows of $\xi_1,...,\xi_s$. The closure $\overline{G}$ of $G$ in $Isom(M,g)$ is a finitely dimensional compact Lie group since it is a closed subgroup of $Isom(M,g)$ which itself is well known to be finitely dimensional and compact. This allows us to define a left invers $a:\Omega^{\bullet}(M)\to\Omega_{\mathfrak{g}}^{\bullet}(M)$ to the inclusion $i:\Omega_{\mathfrak{g}}^{\bullet}(M)\to\Omega^{\bullet}(M)$ by the formula:
$$a(\eta_{i_1}...\eta_{i_k}\alpha):=\eta_{i_1}...\eta_{i_k}\tilde{a}(\alpha),$$
where $\tilde{a}(\alpha)$ denotes the averaging of the form $\alpha$ over the action of $\overline{G}$. We denote the kernel of $a$ by $\Omega_0^{\bullet}(M)$.
\begin{prop}\label{splat} Under the above notation $a$ is a chain map inducing a splitting of the chain complex:
$$\Omega^{\bullet}(M)=\Omega_{\mathfrak{g}}^{\bullet}(M)\oplus\Omega_0^{\bullet}(M),$$
and consequently in cohomology:
$$H_{dR}^{\bullet}(M)=H_{\mathfrak{g}}^{\bullet}(M)\oplus H_0^{\bullet}(M).$$ 
\begin{proof} Let us start by noting that this is in fact well defined as the averaging of a transverse form remains transverse (since the action of $\overline{G}$ preserves $T\mathcal{F}$ due to the fact that the action of $G$ preserves $T\mathcal{F}$) while being basic follows by invariance of the resulting form with respect to the action of $\overline{G}$ and hence $G$.
\newline\indent Secondly, let us prove that this is in fact a cochain map using Proposition \ref{split}:
\begin{eqnarray*} da(\eta_{i_1}...\eta_{i_k}\alpha)&=&\sum_{l=1}^k(-1)^{l-1}\eta_{i_1}...\hat{\eta}_{i_l}...\eta_{i_k}d\eta_{i_l}\tilde{a}(\alpha)+(-1)^k\eta_{i_1}...\eta_{i_k}d\tilde{a}(\alpha)\\
&=&\sum_{l=1}^k(-1)^{l-1}\eta_{i_1}...\hat{\eta}_{i_l}...\eta_{i_k}(d^{(0,1)}\eta_{i_l}+d^{(2,-1)}\eta_{i_l})\tilde{a}(\alpha)+(-1)^k\eta_{i_1}...\eta_{i_k}d\tilde{a}(\alpha)\\
&=&\sum_{l=1}^k(-1)^{l-1}\eta_{i_1}...\hat{\eta}_{i_l}...\eta_{i_k}(d^{(0,1)}\eta_{i_l}\tilde{a}(\alpha)+\tilde{a}((d^{(2,-1)}\eta_{i_l})\alpha))+(-1)^k\eta_{i_1}...\eta_{i_k}\tilde{a}(d\alpha)\\
&=& a(d(\eta_{i_1}...\eta_{i_k}\alpha)).
\end{eqnarray*}
\newline\indent The splitting of chain complexes now follows from the fact that $a$ is a left inverse to $i$.
\end{proof}
\end{prop}
The previous Proposition implies that the spectral sequence $E^{p,q}_r$ converges to $H_{dR}^{\bullet}(M)$ if and only if $H^{\bullet}_{0}(M)=0$. Moreover, this splitting behaves well on $\overline{G}$-invariant forms:
\begin{prop}Under the above notation the restriction of $a$ to $\overline{G}$-invariant forms $\Omega^{\bullet}_{\overline{G}}(M)$ is a chain map inducing a splitting of the chain complex:
$$\Omega_{\overline{G}}^{\bullet}(M)=\Omega_{\mathfrak{g},\overline{G}}^{\bullet}(M)\oplus\Omega_{0,\overline{G}}^{\bullet}(M),$$
and consequently in cohomology
$$H_{dR}^{\bullet}(M)\cong H_{\overline{G}}^{\bullet}(M)=H_{\mathfrak{g},\overline{G}}^{\bullet}(M)\oplus H_{0,\overline{G}}^{\bullet}(M).$$
In particular:
$$H_{\mathfrak{g},\overline{G}}^{\bullet}(M)\cong H_{\mathfrak{g}}^{\bullet}(M),\quad H_{0,\overline{G}}^{\bullet}(M)\cong H_{0}^{\bullet}(M).$$
\begin{proof} We start by noting that for each $i\in\{1,...,s\}$ the operator $\mathcal{L}_{\xi_i}$ preserves the decomposition $\Omega^{\bullet}(M)=\Omega_{\mathfrak{g}}^{\bullet}(M)\oplus\Omega_0^{\bullet}(M)$. For forms in both $\Omega_{\mathfrak{g}}^{\bullet}(M)$ and $\Omega_0^{\bullet}(M)$ this is a simple observation due to the formula $\mathcal{L}_{\xi_i}=d\iota_{\xi_i}+\iota_{\xi_i}d$ and the fact that both $d$ and $\iota_{\xi_i}$ preserve such forms. This implies that a form is $\overline{G}$-invariant if and only if both its components with respect to the splitting $\Omega_{\mathfrak{g}}^{\bullet}(M)\oplus\Omega_0^{\bullet}(M)$ are $\overline{G}$-invariant.
\newline\indent From the previous paragraph we can infer that $a$ takes $\overline{G}$-invariant forms to $\overline{G}$-invariant forms since $a(\alpha_{\mathfrak{g}}+\alpha_0)=\alpha_{\mathfrak{g}}$ which is invariant whenever $\alpha_{\mathfrak{g}}+\alpha_0$ is invariant. From this and the fact that $a$ is a left inverse to $i$ (as chain maps) we have the splittings of $\overline{G}$-invariant forms and cohomology stated in the Proposition. Finally, the isomorphisms:
$$H_{\mathfrak{g},\overline{G}}^{\bullet}(M)\cong H_{\mathfrak{g}}^{\bullet}(M),\quad H_{0,\overline{G}}^{\bullet}(M)\cong H_{0}^{\bullet}(M),$$
are simply the morphisms induced by inclusions $\Omega_{\mathfrak{g},\overline{G}}^{\bullet}(M)\subset \Omega_{\mathfrak{g}}^{\bullet}(M)$ and $\Omega_{0,\overline{G}}^{\bullet}(M)\subset\Omega_{0}^{\bullet}(M)$. It is clear that these are isomorphisms since it is well known that the inclusion of $\overline{G}$-invariant forms induces an isomorphism in cohomology $H_{dR}^{\bullet}(M)\cong H_{\overline{G}}^{\bullet}(M)$ for compact connected Lie groups.
\end{proof}
\end{prop}
\section{Vanishing of $H^{\bullet}_{0}(M)$}\label{Vanish}
We finish the proof of Theorem \ref{MainC} (and consequently Theorem \ref{MainA}) by proving that under its assumptions the group $H^{\bullet}_{0}$ is trivial. In order to do this we note that $H_0^{\bullet}(U)$ is well defined for any saturated open subset $U$ of $M$. Moreover, the splittings defined in the previous sections also work for such a set $U$. We will use this in order to compute $H_0^{\bullet}(M)$ applying the following version of Mayer-Vietoris Theorem:
\begin{lem}\label{MV0} Let $(M,g)$ be a compact Riemannian manifold with a free Killing action of $\mathfrak{g}$. Moreover, let $U$ and $V$ be open saturated subsets. Then the following short sequences of chain complexes is exact:
$$0\to \Omega_0^{\bullet}(U\cup V)\to\Omega_0^{\bullet}(U)\oplus\Omega_0^{\bullet}(V)\to\Omega_0^{\bullet}(U\cap V)\to 0.$$
\begin{proof} The only non-tirival part is the surjectivity of the last non-trivial chain map. This is proved via the same argument as in the proof of the Mayer-Vietoris Theorem in \cite{ttt} after noting that we can construct a partition of unity subordinate to $\{U,V\}$ consisting of basic functions (take any partition of unity suboordinate to $\{U,V\}$ and average it over $\overline{G}$).
\end{proof}
\end{lem}
To apply this lemma effectively we first need to conduct some local computation. We start with the following lemma:
\begin{lem}\label{D2} Under the above notation, $\overline{G}=G\cdot \overline{Z(G)}$ (where the latter group is the closure of the center of $G$ in $\overline{G}$). Moreover, (with respect to any bi-invariant metric on $\overline{G}$) $\mathfrak{g}^{\perp}$ is a Lie subalgebra of $\overline{\mathfrak{g}}$ (Lie algebra of $\overline{G}$) corresponding to a group $A\subset\overline{Z(G)}$.
\begin{proof} We need to prove that the inner automorphisms $I(G)\subset Aut(G)$ form a closed set in order to use Theorem \ref{Dense}. Let $\mathfrak{g}^{\perp}$ be the subspace perpendicular to $\mathfrak{g}$ in $\overline{\mathfrak{g}}$ with respect to a bi-invariant metric on $\overline{G}$. Then a left invariant basis $\{\tilde{\xi}_1,...,\tilde{\xi}_s\}$ of $\mathfrak{g}$ can be extended to a basis of $\overline{\mathfrak{g}}$ using left invariant vector fields $\{\tilde{\mu}_1,...,\tilde{\mu}_l\}$ in $\mathfrak{g}^{\perp}$. Moreover, such vector fields commute with each $\overline{\xi}_i$, since due to the fact that the flows of these vector fields preserve the decomposition $\mathfrak{g}^{\perp}\oplus \mathfrak{g}$, we have $[\tilde{\mu}_i,\tilde{\xi_j}]\in \mathfrak{g}\cap \mathfrak{g}^{\perp}=\{0\}$ and $[\tilde{\mu}_i,\tilde{\mu_j}]\in \mathfrak{g}^{\perp}$. From this, we can write every element of $\overline{G}$ as $ag$ with $g\in G$ and $a\in A$ such that $\overline{G}=G\cdot A$ and such that all elements in $A$ commute with elements in $G$. We can conclude from this that $I(\mathfrak{g})\cong I(G)$ is an embedded subgroup of $I(\overline{\mathfrak{g}})\cong I(\overline{G})$ (since any element $\xi\in\mathfrak{g}$ acts trivially on $\mathfrak{g}^{\perp}$ and hence $e^{ad_{\xi}}|_{\mathfrak{g}^{\perp}}=Id_{\mathfrak{g}^{\perp}}$). Moreover, the embedding $E:I(G)\to I(\overline{G})$ is in fact bijective (since conjugation by $ag$ has to agree with conjugation by $g$ on $G$ and consequently the same is true on $\overline{G}$). Hence, $I(G)$ is isomorphic as Lie groups with $I(\overline{G})$ which implies that $I(G)$ is compact (in particular, closed in $Aut(G)$). Hence, by Theorem \ref{Dense} we have $\overline{G}=G\cdot \overline{Z(G)}$. In particular, $A\subset\overline{Z(G)}$ since it commutes with all the elements of $G$.
\end{proof}
\end{lem}
This enables us to provide the following result:
\begin{prop}\label{loc} Let $(M,g)$ be a compact Riemannian manifold with a free Killing action of $\mathfrak{g}$ and let $V$ be a saturated open subset admitting a slice through some orbit contained in $V$. Then $H_0(V)=0$.
\begin{proof}
It is sufficient to show that each class in $H^*_{dR}(V)$ has a representative in $\Omega_{\mathfrak{g},\overline{G}}(V)$ due to Proposition \ref{splat}. Since $V$ is a neighbourhood admitting a slice through some point $x\in V$ its cohomology is the same as that of $\overline{G}x\cong\overline{G}\slash \overline{G}_x$. Moreover, a form in $\Omega_{\mathfrak{g},\overline{G}}(\overline{G}\slash\overline{G}_x)$ can be extended to a form in $\Omega_{\mathfrak{g},\overline{G}}(V)$ using the Slice Theorem so that it represents the corresponding cohomology class.
\newline\indent Let $T\supset G$ denote a connected Lie group (with Lie algebra $\mathfrak{t}$) which forms a transversal for the foliation induced by the action of $\overline{G}_x$ on $\overline{G}$ (with Lie algebras respectively $\overline{\mathfrak{g}}_x$ and $\overline{\mathfrak{g}}$) formed by taking the Lie group corresponding to $\mathfrak{g}\oplus span(\tilde{\mu}_{i_1},...,\tilde{\mu}_{i_r})$ for an appropriate subset $\{\tilde{\mu}_{i_1},...,\tilde{\mu}_{i_r}\}$ of a basis of $\mathfrak{g}^{\perp}$. We consider $T$ with its Lie group topology (the subgroup topology would also work for the argument). We shall show that every cohomology class in $H^{\bullet}_{dR}(\overline{G}\slash\overline{G}_x)$ has a representative that lifts to a $T$-bi-invariant form on $T$. Firstly, let us note that $(\overline{G}_x)_0$ is abelian since if there where any nonvanishing bracket of elements from $\overline{\mathfrak{g}}_x$ it would have to be contained (by Lemma \ref{D2}) in $\mathfrak{g}\subset\mathfrak{t}$ which gives a contradiction since $\mathfrak{t}$ and $\overline{\mathfrak{g}}_x$ meet transversaly (and are of complementary dimensions). We put, $\Gamma_T:=T\cap (\overline{G}_x)_0\cap\overline{Z(G)}\subset T$ and show that $T\slash\Gamma_T$ is a compact Lie group. To this end we first need to show that $(\overline{G}_x)_0\cap T\slash\Gamma_T$ is finite. To see this assume to the contrary this is not the case and consider the compact group $T\slash (T\cap\overline{Z(G)})\cong \overline{G}\slash\overline{Z(G)}$ (the natural mapping between the two is well known to be a continuous bijection and computing the differential at identity proves that it is in fact an isomorphism of Lie groups). Consequently, the space:
$$T\slash (T\cap\overline{Z(G)})\bigg/ ((\overline{G}_x)_0\cap T)\slash\Gamma_T\cong \overline{G}\slash\overline{Z(G)}\bigg/ (\overline{G}_x)_0\slash ((\overline{G}_x)_0\cap\overline{Z(G)}),$$
has to be a manifold (due to the right side) despite being covered by a compact manifold with infinite fiber (left side). Hence, we get our desired contradiction.
\newline\indent Since, $T$ covers $\overline{G}\slash (\overline{G}_x)_0$ and the group $(\overline{G}_x)_0\cap T\slash\Gamma_T$ is finite, we conclude that $T\slash\Gamma_T$ is a finite cover of:
$$\overline{G}\slash (\overline{G}_x)_0\cong T\slash ((\overline{G_x})_0\cap T)\cong (T\slash\Gamma_T)\slash ((\overline{G}_x)_0\cap T\slash\Gamma_T),$$ 
where again we can conclude that the natural mappings between the above spaces are diffeomorphism since they are bijective continuous and the differential at identity (and hence at any other point by translation by elements of $T$ is an isomorphism). Hence, due to the compactness of the later space $T\slash\Gamma_T$ is compact itself. Hence: 
$$T\slash\Gamma_T\to T\slash ((\overline{G_x})_0\cap T)\to T\slash (\overline{G_x}\cap T)\cong\overline{G}\slash\overline{G}_x,$$
forms a tower of finite coverings. Using that finite coverings induce injections in de Rham cohomology as well as the fact that $T\slash\Gamma_T$ is a compact connected Lie group with Lie algebra $\mathfrak{t}$ we deduce that to each cohomology class on  $\overline{G}\slash\overline{G}_x$ we can associate a unique $T$-bi-invariant form on $T$. Finally, let us note that each such form descends to a form on $\overline{G}\slash\overline{G}_x\cong T\slash (\overline{G}_x\cap T)$ representing the initial cohomology class.
\newline\indent Using the previous paragraph, let $\alpha$ be a $T$-bi-invariant closed $k$-form on $T$. Let $\{X_1,...,X_s\}$ be a basis of left invariant vector fields tangent to $G\subset T$ coinciding at identity with the basis $\{\overline{\xi}_1,...,\overline{\xi}_s\}$ of right invariant vector fields. Put $\overline{\eta}_i:=\pi^*\eta_i|_{\overline{G}x}=(\pi^*g_T)(\overline{\xi}_i,\bullet)$. Note that $\mathcal{L}_{X_j}\overline{\eta}_i=0$ since $\pi^*g_T$ is bi-invariant. This gives:
$$0=\mathcal{L}_{X_j}\alpha=\mathcal{L}_{X_j}\alpha_0+\sum_{1\leq i_1<...<i_l\leq s} \overline{\eta}_{i_1}....\overline{\eta}_{i_l}\mathcal{L}_{X_j}\alpha_{i_1,...,i_l}.$$
But since $G$ is a dense subgroup (in $\overline{G}\supset T$) the bundle $TG$ tangent to the layers of $G$ is preserved by multiplication on both sides and consequently $\mathcal{L}_{X_j}\alpha_0$ and $\mathcal{L}_{X_j}\alpha_{i_1,...,i_l}$ are again normal. This means that the above equality holds only if $\mathcal{L}_{X_j}\alpha_0=0$ and $\mathcal{L}_{X_j}\alpha_{i_1,...,i_l}=0$. Since this holds for all $X_j$ the forms $\alpha_0$ and $\alpha_{i_1,...,i_l}$ are all basic (with respect to the orbits of the left $G$-action). Consequently, if $\alpha$ represents a pullback of a form on $\overline{G}\slash\overline{G}_x$ then this form belongs to $\Omega_{\mathfrak{g},\overline{G}}(\overline{G}x)$. By extending this form to $V$ the proof is finished.
\end{proof}
\end{prop}
With this we are able to prove the following Proposition which finishes the proof of Theorems \ref{MainC} and \ref{MainA}.
\begin{prop} For any compact manifold $(M^{n+s},g)$ with a Killing free action of an $s$-dimensional Lie algebra $\mathfrak{g}$, the vector spaces $H_0^{\bullet}(M)$ are all trivial.
\begin{proof} We are going to proceed with the proof by induction on the number of strata. More precisely, we prove the following statement by induction:
\newline\newline\indent\emph{Statement}. Let $N^{n+s}$ be a saturated open submanifold of some compact manifold $M^{n+s}$ with a free Killing $\mathfrak{g}$-action. Moreover, let us assume that $N$ is $\overline{G}$-equivariantly (where the closure is taken in $Isom(M,g)$) diffeomorphic to an interior of a submanifold with boundary $\tilde{N}$ of $M$. If, the action of $\overline{G}$ on $N$ has at most $k$ distinct strata then $H^{\bullet}_{0}(N)=0$.
\newline\newline\indent For $k=1$ we can use the Slice Theorem to see that $N\slash\overline{G}$ is a manifold. We can take $\mathcal{U}$ to be the cover of $N$ consisting of preimages of sets from a finite good open cover of $N\slash\overline{G}$. This implies that $\mathcal{U}$ is an open saturated cover of $N$ such that any intersection of sets from that cover admits a slice. We can now use Lemma \ref{MV0} and Proposition \ref{loc} to conclude that $H^{\bullet}_{0}(N)=0$ (by induction on the cardinality of the cover).
\newline\indent For the induction step let us assume that the statement is true for all positive integers up to $k$ and let us prove the statement for $k+1$. Let $K_{k+1}$ be the final strata of $N$ (in particular the orbits in it are not quotients of orbits from any other orbit type). Moreover, let us take some equivariant tubular neighbourhood $U$ of $K_{k+1}$. It can be seen from the induction hypothesis that $H^{\bullet}_0(K_{k+1})=0$ since, one can see using the slice Theorem that $\overline{K}_{k+1}$ is an immersed (since it might have intersections) compact submanifold in $M$ and consequently $K_{k+1}$ satisfies the assumptions of the induction hypothesis (where the place of $M$ is taken by a compact manifold which is the source of the immersion with image $\overline{K}_{k+1}$). This implies that, $H^{\bullet}_0(U)=0$ (since otherwise a restriction of the non-trivial class from $U$ would be a nontrivial class on $K_{k+1}$).
\newline\indent We wish to now prove the same for $N\backslash K_{k+1}$. To this end note that by induction hypothesis $H^{\bullet}_0 (N\backslash \overline{U})=0$. On the other hand, $N\backslash \overline{U}$ is equivariantly diffeomorphic by Proposition \ref{tub} to $N\backslash K_{k+1}$ and consequently by applying Lemma \ref{MV0} to $U$ and $N\backslash K_{k+1}$ we get $H^{\bullet}_0(N)=0$ (noting that $H^{\bullet}_0(U\cap (N\backslash K_{k+1}))=0$, since $U\cap (N\backslash K_{k+1})$ satsfies the assumptions of the induction hypothesis).
\end{proof}
\end{prop}
\section{Applications to low codimension actions}
In this section we provide some sample applications for our spectral sequence for low codimension actions. These results for the most part can be thought of as a generalization of the Wang long exact sequence. We start however with a simple consequence of the spectral sequence in general (which is a genralization of a similar result from \cite{My3}).
\begin{prop}\label{Hor}  Let $(M^{n+s},g)$ be a compact orientable manifold with a Killing free action of an $s$-dimensional Lie algebra $\mathfrak{g}$. Then the foliation induced by the action is homologically orientable.
\begin{proof} It is well known (cf. \cite{E1}) that the top basic cohomology of a Riemannian foliation on a compact manifold is either $0$ or $\mathbb{R}$. In this case it cannot be $0$ since then we could compute from our spectral sequence that $H^{n+s}_{dR}(M)\cong E^{n,s}_2=0$ which is a contradiction with the orientability of $M$. Hence, the top basic cohomology is isomorphic to $\mathbb{R}$ which means that the foliation is homologically orientable.
\end{proof}
\end{prop}
\begin{tw} Let $(M^{1+s},g)$ be a compact connected oriented Riemannian manifold with a free Killing action of an $s$-dimensional Lie algebra $\mathfrak{g}$. Then:
$$H^{\bullet}_{dR}(M)\cong H^{\bullet}(\mathfrak{g})\oplus H^{\bullet-1}(\mathfrak{g}).$$
In particular, no such action exists when $M$ is simply connected (or even when $H^1_{dR}(M)=0$).
\begin{proof} Since $E^{0,q}_2$ and $E^{1,q}_2$ are the only non-zero columns on the second page we conclude that the sequence degenerates at this page. Due to Proposition \ref{Hor} we can use the basic version of Poincar\'e duality to conclude from $E_2^{p,q}\cong E^{p,q}_{\infty}$ that the cohomology of $M$ is as prescribed. Moreover, $H^{0}(\mathfrak{g})\neq 0$ together with the above splitting implies that $H^{1}_{dR}(M)\neq 0$. \end{proof}
\end{tw}
\begin{tw} Let $(M^{2+s},g)$ be a compact connected oriented Riemannian manifold with a free Killing action of an $s$-dimensional Lie algebra $\mathfrak{g}$ and $H^{1}_{dR}(M)=0$. Then we have a long exact sequence:
$$...\to H^{k}_{dR}(M)\to H^{k}(\mathfrak{g})\to H^{k-1}(\mathfrak{g})\to  H^{k+1}_{dR}(M)\to ...$$.
\begin{proof} We use the well known fact that the morphism $i^*:H^1(M\slash\mathcal{F})\to H^{1}_{dR}(M)$ induced by the inclusion of forms is a monomorphism to conclude that $H^1(M\slash\mathcal{F})=0$ and hence only the $0$-th and the $2$-nd column of $E_2^{p,q}$ are nonzero. Moreover, due to Proposition \ref{Hor} we can use the basic version of Poincar\'e duality to conclude that  $E^{0,q}_2\cong H^{q}(\mathfrak{g})$ and $E^{2,q}_2\cong H^{q}(\mathfrak{g})$. Hence, we can use the derivative on the second page to arrive at the above long exact sequence.
\end{proof}
\end{tw}

\begin{tw} Let $(M^{3+s},g)$ be a compact connected oriented Riemannian manifold with a free Killing action of an $s$-dimensional Lie algebra $\mathfrak{g}$ and $H^{1}_{dR}(M)=0$. Then we have a long exact sequence:
$$...\to H^{k}_{dR}(M)\to H^{k}(\mathfrak{g})\to H^{k-2}(\mathfrak{g})\to  H^{k+1}_{dR}(M)\to ...$$
Moreover, $0\cong H^{1}(\mathfrak{g})\cong\mathfrak{g}\slash [\mathfrak{g},\mathfrak{g}]$.
\begin{proof} We again use the well known fact that the morphism $i^*:H^1(M\slash\mathcal{F})\to H^{1}_{dR}(M)$ induced by the inclusion of forms is a monomorphism to conclude that $H^1(M\slash\mathcal{F})=0$. Moreover, we can use Proposition \ref{Hor} to conclude from the basic version of Poincar\'e duality that $H^2(M\slash\mathcal{F})=0$ and hence, that  $E^{0,q}_2\cong H^{q}(\mathfrak{g})$ and $E^{3,q}_2\cong H^{q}(\mathfrak{g})$ are the only non-zero columns. Hence, we can use the derivative on the third page to arrive at the above long exact sequence. The final assertion follows from the fact that there are no non-zero coboundary operators with source or target on $E^{0,1}_r$ for $r\geq 2$ and consequently $H^1(\mathfrak{g})\cong E^{0,1}_2\cong E^{0,1}_{\infty}\cong H^1_{dR}(M)=0$.
\end{proof}
\end{tw}
\begin{rem} The assumption $0\cong \mathfrak{g}\slash [\mathfrak{g},\mathfrak{g}]$ is pretty restrictive and prevents a lot of groups from acting locally freely and isometricaly with codimension $3$ on compact orientable manifolds with trivial first cohomology. Among these it is worth mentioning that $\mathbb{R}^s$ cannot act in such a way on such manifolds. This can be viewed as a generalization of the fact that there are no nowhere vanishing vector fields on simply connected compact orientable $4$-manifolds (which follows immediately from the Euler characteristic).
\end{rem}
It is also worth noting that some similar results can be achieved from the study of our spectral sequence for actions of higher codimension if we make additional assumptions on the vanishing of low degree cohomology of $M$ and Lie algebra cohomology of $\mathfrak{g}$.

\bigskip
\Small{\textbf{Statements and Declarations} The authors declare no competing interests.}
\bigskip

\Small{\textbf{Availability of data and material} Data sharing not applicable to this article as no datasets were generated or analysed during the current study.}

\end{document}